\documentclass[letterpaper]{amsart}
\usepackage[english]{babel}
\usepackage{amsmath}
\usepackage{amssymb}
\usepackage{amsfonts}
\usepackage{amsthm}
\usepackage{mathtools}
\usepackage{pstricks}
\usepackage{pst-node}
\usepackage{pst-plot}
\usepackage{pst-3dplot}
\usepackage{booktabs}
\usepackage{multirow}
\usepackage{subfigure}
\usepackage{pdflscape}
\usepackage{url}
\usepackage[all,cmtip]{xy}
\allowdisplaybreaks[3]
\newcommand{\ds}{\displaystyle}
\renewcommand{\subset}{\subseteq}
\renewcommand{\supset}{\supseteq}

\newcommand{\R}{\ensuremath{\mathbf{R}}}

\newcommand{\Rplus}{\ensuremath{\mathbf{R}_{\geqslant0}}}

\newcommand{\Rplusnul}{\ensuremath{\mathbf{R}_{>0}}}

\newcommand{\Rn}{\ensuremath{\mathbf{R}^n}}
\newcommand{\Rplusn}{\ensuremath{\mathbf{R}_{\geqslant0}^n}}

\newcommand{\N}{\ensuremath{\mathbf{N}}}
\newcommand{\Nnul}{\ensuremath{\mathbf{N}\setminus\{0\}}}
\newcommand{\Nn}{\ensuremath{\mathbf{N}^n}}
\newcommand{\Z}{\ensuremath{\mathbf{Z}}}
\newcommand{\Zn}{\ensuremath{\mathbf{Z}^n}}

\newcommand{\Zplus}{\ensuremath{\mathbf{Z}_{\geqslant0}}}

\newcommand{\Zplusn}{\ensuremath{\mathbf{Z}_{\geqslant0}^n}}

\newcommand{\Zplusnul}{\ensuremath{\mathbf{Z}_{>0}}}

\newcommand{\Q}{\ensuremath{\mathbf{Q}}}

\newcommand{\C}{\ensuremath{\mathbf{C}}}

\newcommand{\Ccrossn}{\ensuremath{(\mathbf{C}^{\times})^n}}
\newcommand{\Qp}{\ensuremath{\mathbf{Q}_p}}

\newcommand{\Zp}{\ensuremath{\mathbf{Z}_p}}

\newcommand{\Fp}{\ensuremath{\mathbf{F}_p}}

\newcommand{\A}{\ensuremath{\mathbf{A}}}

\renewcommand{\N}{\Zplus}
\renewcommand{\Nnul}{\Zplusnul}
\renewcommand{\Nn}{\Zplusn}
\renewcommand{\phi}{\varphi}

\DeclareMathOperator{\aff}{aff}

\DeclareMathOperator{\supp}{supp}
\DeclareMathOperator{\mult}{mult}

\DeclareMathOperator{\odiv}{div}
\DeclareMathOperator{\codim}{codim}
\DeclareMathOperator{\Vol}{Vol}
\DeclareMathOperator{\res}{res}

\newcommand{\I}{\ensuremath{\mathcal{I}}}
\newcommand{\Ig}{\ensuremath{(\mathcal{I},gdx)}}

\newcommand{\ZIg}{\ensuremath{Z_{\mathcal{I},gdx}}}

\newcommand{\Ztopf}{\ensuremath{Z^{\mathrm{top}}_f}}

\newcommand{\ZtopI}{\ensuremath{Z^{\mathrm{top}}_{\mathcal{I}}}}

\newcommand{\ZtopIg}{\ensuremath{Z^{\mathrm{top}}_{\mathcal{I},gdx}}}
\newcommand{\ZtopoIg}{\ensuremath{Z^{\mathrm{top},0}_{\mathcal{I},gdx}}}
\newcommand{\ZtopglIg}{\ensuremath{Z^{\mathrm{top,gl}}_{\mathcal{I},gdx}}}

\newcommand{\ZtopIgs}{\ensuremath{Z^{\mathrm{top}}_{\mathcal{I},gdx}(s)}}
\newcommand{\ZtopglIgs}{\ensuremath{Z^{\mathrm{top,gl}}_{\mathcal{I},gdx}(s)}}
\newcommand{\ZtopoIgs}{\ensuremath{Z^{\mathrm{top},0}_{\mathcal{I},gdx}(s)}}

\newcommand{\mI}{\ensuremath{m_{\mathcal{I}}}}
\newcommand{\bI}{\ensuremath{b_{\mathcal{I}}}}
\newcommand{\bIs}{\ensuremath{b_{\mathcal{I}}(s)}}

\newcommand{\Gg}{\ensuremath{\Gamma_g}}
\newcommand{\GI}{\ensuremath{\Gamma_{\mathcal{I}}}}

\newcommand{\Gglg}{\ensuremath{\Gamma^{\mathrm{gl}}_g}}

\newcommand{\Dg}{\ensuremath{\Delta_g}}
\newcommand{\DI}{\ensuremath{\Delta_{\mathcal{I}}}}

\newcommand{\DIg}{\ensuremath{\Delta_{\mathcal{I},g}}}

\newcommand{\DIt}{\ensuremath{\Delta_{\mathcal{I}}(\tau)}}

\theoremstyle{plain}
\newtheorem{theorem}{Theorem}[section]

\newtheorem{proposition}[theorem]{Proposition}

\newtheorem{question}[theorem]{Question}
\theoremstyle{definition}
\newtheorem{definition}[theorem]{Definition}

\theoremstyle{remark}
\newtheorem{remark}[theorem]{Remark}

\title[Zeta functions and $b$-functions for ideals in dimension two]{Zeta functions and Bernstein--Sato polynomials\\for ideals in dimension two}
\author{Bart Bories}
\address{Department of Mathematics, KU Leuven, Celestijnenlaan 200b -- box 2400, 3001 Leuven, Belgium}
\email{bart.bories@wis.kuleuven.be}
\date{\today}

\begin{document}
\begin{abstract}
For a nonzero ideal $\I\lhd\C[x_1,\ldots,x_n]$, with $0\in\supp\I$, a generalization of a conjecture of Igusa--Denef--Loeser predicts that every pole of its topological zeta function is a root of its Bernstein--Sato polynomial. However, typically only a few roots are obtained this way. Following ideas of Veys \cite{Veys07}, we study the following question. Is it possible to find a collection $\mathcal{G}$ of polynomials $g\in\C[x_1,\ldots,x_n]$, such that, for all $g\in\mathcal{G}$, every pole of the topological zeta function associated to \I\ and the volume form $gdx_1\wedge\cdots\wedge dx_n$ on the affine $n$-space, is a root of the Bernstein--Sato polynomial of \I, and such that all roots are realized in this way. We obtain a negative answer to this question, providing counterexamples for monomial and principal ideals in dimension two, and give a partial positive result as well.
\end{abstract}

\maketitle

\setcounter{tocdepth}{1}
\tableofcontents

\section*{Introduction}
The topological, motivic, and $p$-adic Igusa zeta function are invariants originally associated to polynomials or analytic functions in several variables over \C, an arbitrary field of characteristic zero, and a $p$-adic field, respectively. Intriguing conjectures, motivated by results in classical complex integration theory, link the poles of these rational functions to monodromy eigenvalues and the roots of Bernstein--Sato polynomials.

For an analytic function germ $f:(\C^n,0)\to(\C,0)$, the \lq Monodromy Conjecture\rq\ of Denef and Loeser predicts that every pole of the topological zeta function \Ztopf\ of $f$ induces one of its local monodromy eigenvalues. Another conjecture of Denef and Loeser asserts that every such pole is a root of the Bernstein--Sato polynomial (also called $b$-function) $b_f$ of $f$. (There are motivic and $p$-adic versions of these conjectures, the $p$-adic ones due to Igusa.)

These three zeta functions (topological, motivic, and $p$-adic) have been generalized in a straightforward way to polynomial mappings or ideals in polynomial rings. Not so obvious was the generalization of the concepts of local monodromy and \mbox{$b$-functions} to several polynomials, this was done by Verdier \cite{Verdier} for local monodromy and by Sabbah \cite{Sa87} and more recently by Budur et al.\ \cite{BMS06bis} for Bernstein--Sato polynomials. The conjectures we mentioned can still be stated in this broader context.

The Monodromy Conjecture has been proved in full generality for $n=2$ by Loeser \cite{Loe88} (for one analytic function, originally in the context of $p$-adic Igusa zeta functions) and by Van Proeyen and Veys \cite{VVmcid2} (for polynomial ideals); in higher dimension there are various partial results, e.g., \cite{ACLM02,ACLM05,BMTmcha,Kimura,LVmcndss,LV09,Loe90,RodrVeys,Vey93,Vey06}. The converse of the conjecture is certainly not true. In fact, generally, only a few monodromy eigenvalues are induced by poles of the corresponding zeta function. This led Veys to consider zeta functions associated to an analytic function germ $f:(\C^n,0)\to(\C,0)$ and a differential $n$-form $\omega$, and pose the following question. Can one find a family $\Omega$ of differential forms, such that, for all $\omega\in\Omega$, every pole of the topological zeta function $Z^{\mathrm{top}}_{f,\omega}$ of $f$ and $\omega$ induces a monodromy eigenvalue of $f$, and such that all monodromy eigenvalues of $f$ are obtained this way? In \cite{NV10bis} N{\'e}methi and Veys give an affirmative answer to this question for an arbitrary curve singularity $f$.

The situation is very similar for the stronger\footnote{It is well-known that roots of the Bernstein--Sato polynomial induce local monodromy eigenvalues \cite{malgrange}.} conjecture of Igusa--Denef--Loeser, expecting zeta function poles to be roots of the corresponding $b$-function. The conjecture has been verified by Loeser for one analytic function in two variables \cite{Loe88} and (under extra technical conditions) for non-degenerated polynomials in several variables \cite{Loe90}, but for polynomial ideals, not much is known yet. In \cite{HMY07} Howald, Musta{\c{t}}{\u{a}}, and Yuen proved the conjecture for monomial ideals in arbitrary dimension in the context of $p$-adic Igusa zeta functions, but their argument can easily be adapted to cover the topological and motivic versions of the conjecture as well (see Remark~\ref{conj2monidtop}). Again, examples show that the poles of a zeta function generally cover only a small part of the roots of the corresponding Bernstein--Sato polynomial. The positive result of Veys and N{\'e}methi for the Monodromy Conjecture provides motivation to study the analogous question for the second conjecture. In this paper, we conclude negatively for the topological zeta function of a monomial ideal $\I\lhd\C[x_1,x_2]$ in dimension two, and provide a counterexample for the one polynomial case as well. A partial positive result is obtained in Section~\ref{positiveresults}.

\section*{Acknowledgments}
The author would like to thank Wim Veys for proposing the problem and for many useful suggestions and comments about the paper.

\section{Preliminaries}
\subsection{The topological zeta function of an ideal and a polynomial volume form}
\begin{definition}[Embedded resolution]\label{def_er}
Let \I\ be an ideal of $\C[x]=\C[x_1,\ldots,x_n]$, with $0=(0,\ldots,0)\in\supp\I$, and $g$ a nonzero polynomial in $\C[x]$. In $\A^n(\C)$ we consider the support of \I\ and the volume form $gdx=gdx_1\wedge\cdots\wedge dx_n$. By an embedded resolution of \Ig, we mean in this paper a composition
\begin{equation*}
h=h_1\circ\cdots\circ h_t:Y=Y_t\to Y_0=\A^n(\C)
\end{equation*}
of blowing-ups $h_i:Y_i\to Y_{i-1}$ in smooth centra $C_{i-1}\subset Y_{i-1}$, such that
\begin{enumerate}
\item $h$ is a principalization of \I; i.e.,
\begin{enumerate}
\item for all $i\in\{1,\ldots,t\}$, the exceptional divisor of the morphism $h_1\circ\cdots\circ h_i:Y_i\to Y_0$ is a simple normal crossings divisor on $Y_i$, having moreover simple normal crossings with $C_i$ for $i\in\{1,\ldots,t-1\}$, and
\item the total transform (pull-back) $h^{\ast}(\I)$ of \I\ is the ideal of a simple normal crossings divisor $F$ on $Y$;
\end{enumerate}
\item $h$ is an embedded resolution of singularities of the locus $g^{-1}(0)\subset\A^n(\C)$, where we denote by $G=\odiv g\circ h$ the (principal) simple normal crossings divisor on $Y$ associated to $g\circ h$; and
\item $F+G$ is (still) a simple normal crossings divisor on $Y$.
\end{enumerate}
\end{definition}

\begin{remark}
With the notations of Definition~\ref{def_er}, let $h$ be an embedded resolution of \Ig, and denote by $E$ the exceptional divisor of $h$ on $Y$. If the support of \I\ has no components of codimension one, then $F$ is a linear combination (with nonnegative integral coefficients) of the prime divisors of $E$. Otherwise, we can factor $h^{\ast}(\I)$ into a product of two (locally principal) ideals: the support of the first one is contained in the exceptional locus $\supp E$, while the support of the second ideal (called the weak transform of \I) is the union of the strict transforms of the codimension one irreducible components of the support of~\I.

Analogously, the divisor $G$ can be written as a sum $G_{\text{e}}+G_{\text{s}}$, where $\supp G_{\text{e}}\subset\supp E$ and the support of $G_{\text{s}}$ is the strict transform of $g^{-1}(0)$ under $h$.
\end{remark}

\begin{remark}
Let \I\ and $g$ be as above. An embedded resolution of \Ig\ can be obtained as the composition $h'\circ h''$ of a principalization $h'$ of \I\ and an embedded resolution of singularities $h''$ of the union of the supports of $(h')^{\ast}(\I)$ and $(h')^{\ast}(g)=g\circ h'$. The existence of both maps is guaranteed by Hironaka's Theorem.
\end{remark}

\begin{definition}[Numerical data of an embedded resolution]\label{defnumdata}
We use the notations of Definition~\ref{def_er}. Let $h$ be an embedded resolution of \Ig. Denote by $E_i$, $i\in S$, the prime divisors of $F+G$, i.e., the irreducible components of the support of $F+G$, and put $\mathcal{S}=\{E_i\mid i\in S\}$. On $Y$ we consider the ideal sheaf $h^{\ast}(\I)=\mathcal{O}_Y(-F)\lhd\mathcal{O}_Y$ and the pull-back $h^{\ast}(gdx)=h^{\ast}(g)h^{\ast}(dx)$ of the volume form $gdx$. For each $E_i\in\mathcal{S}$, denote by $N_i\geqslant0$ and $\nu_i-1\geqslant0$ the multiplicities of $F$ and $h^{\ast}(gdx)$, respectively, along $E_i$. The $(N_i,\nu_i)$ are called the numerical data of $E_i$. We usually write $E_i(N_i,\nu_i)$ to indicate that $E_i$ has numerical data $(N_i,\nu_i)$.

For all $b\in Y$, we then have the following. Suppose $E_{i_1},\ldots,E_{i_r}\in\mathcal{S}$, with $0\leqslant r\leqslant n$, are the (distinct) irreducible components of $\supp(F+G)$ passing through $b$. Then there exists an open neighborhood $U\subseteq Y$ of $b$ and local coordinates $(y_1,\ldots,y_n)$ on $U$, centered at $b$, such that $y_1=0,\ldots,y_r=0$ are local equations of $E_{i_1},\ldots,E_{i_r}$, respectively, on $U$, and such that
\begin{align*}
h^{\ast}(\I)&=\left(\varepsilon(y)\prod\nolimits_{j=1}^ry_j^{N_{i_j}}\right),\\
F&=\sum_{j=1}^rN_{i_j}E_{i_j},\quad\text{and}\\
h^{\ast}(gdx)&=\eta(y)\prod_{j=1}^ry_j^{\nu_{i_j}-1}dy_1\wedge\cdots\wedge dy_n
\end{align*}
on $U$, for some units $\varepsilon(y)$ and $\eta(y)$ in the local ring $\mathcal{O}_{Y,b}$ of $Y$ at $b$.
\end{definition}

\begin{definition}[Topological zeta function]\label{def_topzetafie}
We use the notations of Definition~\ref{defnumdata}. For $I\subset S$, denote $E_I=\bigcap_{i\in I}E_i$ and $E_I^{\circ}=E_I\setminus\bigcup_{i\not\in I}E_i$.

Based on the embedded resolution $h$ and the corresponding numerical data, one associates to \I\ and $gdx$ the (local) topological zeta function
\begin{equation}\label{def_topzetafie_eq}
\ZtopIg=\ZtopoIg:s\mapsto\ZtopIgs=\sum_{I\subset S}\chi(E_I^{\circ}\cap h^{-1}(0))\prod_{i\in I}\frac{1}{\nu_i+N_is},
\end{equation}
which is a meromorphic function in the complex variable $s$. Here $\chi(\cdot)$ denotes the topological Euler--Poincar\'e characteristic. We obtain a global version \ZtopglIg\ of the topological zeta function replacing $E_I^{\circ}\cap h^{-1}(0)$ by $E_I^{\circ}$ in \eqref{def_topzetafie_eq}.
\end{definition}

\begin{remark}
Since the topological zeta function only depends on the prime divisors that meet $h^{-1}(0)$, to define the function, it is actually sufficient to consider a local embedded resolution of \Ig, i.e., an embedded resolution in some open neighborhood of the origin in $\A^n(\C)$.
\end{remark}

\begin{proposition}
The topological zeta function $\ZtopIg$ does not depend on the choice of embedded resolution of \Ig.
\end{proposition}

\begin{proof}
The topological zeta function was originally defined by Denef and Loeser in \cite{DL92} for a principal ideal $\I=(f)$ and the usual volume form $dx=dx_1\wedge\cdots\wedge dx_n$. Considering it as a kind of limit of Igusa's $p$-adic zeta functions, they proved the remarkable fact that the topological zeta function $\Ztopf=Z^{\mathrm{top}}_{(f),dx}$ does not depend on the embedded resolution (which is just an embedded resolution of singularities of $(f,dx)$ in this case) by which it is defined. Later they obtained the same result considering the topological zeta function as a specialization of the motivic zeta function \cite{DL98}. Veys and Z{\'u}{\~n}iga-Galindo generalized this argument to arbitrary ideals \cite[(2.4)]{VZ08}. Another way to prove it (for ideals) is by comparing two embedded resolutions by means of the Weak Factorization Theorem of Abramovich et al.\ \cite{TWFTref}. One can further adapt these arguments to prove the independence result in our setting of ideals and polynomial volume forms.
\end{proof}

We use the notations of Definition~\ref{def_topzetafie} and we put $\mathcal{S}=\{E_i\mid i\in S\}$. The candidate poles of \ZtopIg\ are the negative rational numbers $-\nu_i/N_i$; $i\in S$, $N_i\neq0$. So each prime divisor $E_i\in\mathcal{S}$, with $N_i\neq0$, gives rise to a candidate pole. In general, some of these numbers will be poles and others will not. We define the expected order of a candidate pole $s_0$ as the maximal number $e$ of prime divisors $E_{i_1},\ldots,E_{i_e}\in\mathcal{S}$, with associated candidate pole $-\nu_{i_j}/N_{i_j}=s_0$; $j=1,\ldots,e$; for which $E_{i_1}\cap\cdots\cap E_{i_e}\cap h^{-1}(0)\neq\emptyset$.
%
%
As the $E_i$ have normal crossings, the expected order of a candidate pole is at most $n$. Clearly, the order as a pole of a candidate pole is always less than or equal to its expected order.

Let us now look closer to the dimension two case. Since for $s_0\in\C$, the limit $\lim_{s\to s_0}(s-s_0)^2\ZtopIgs$ equals the sum of the positive contributions of the intersection points in $h^{-1}(0)\subset Y$ of two prime divisors in $\mathcal{S}$ with associated candidate pole $s_0$, it follows that a candidate pole of expected order two is always a pole of order two.

Let $E_i\in\mathcal{S}$ be a prime divisor with numerical data $(N_i,\nu_i)$, $N_i\neq0$, and an associated candidate pole $s_0=-\nu_i/N_i$ of expected order one. We define the contribution of $E_i$ to the residue of \ZtopIg\ at $s_0$, as the residue at $s_0$ of that part of the sum, in the definition of \ZtopIgs, that depends on $E_i$, i.e., as
\begin{equation*}
\lim_{s\to s_0}(s-s_0)\sum_{I\subset S,\ I\ni i}\chi(E_I^{\circ}\cap h^{-1}(0))\prod_{i\in I}\frac{1}{\nu_i+N_is}.
\end{equation*}
Clearly, the residue of \ZtopIg\ at $s_0$ is the sum of the contributions to this residue of the prime divisors in $\mathcal{S}$ with associated candidate pole $s_0$, and $s_0$ is a pole of \ZtopIg\ if and only if this residue is different from zero. Let $E_{i_1},\ldots,E_{i_r}$ be the prime divisors in $\mathcal{S}\setminus\{E_i\}$ that intersect $E_i$, and put $m_j=\#(E_i\cap E_{i_j})$; $j=1,\ldots,r$. Note that $m_j=1$ if both $E_i$ and $E_{i_j}$ are exceptional. Denote by $m=\sum_j m_j$ the total number of intersection points between $E_i$ and other divisors in $\mathcal{S}$, and put $\alpha_j=\nu_{i_j}+s_0N_{i_j}=\nu_{i_j}-(\nu_i/N_i)N_{i_j}$; $j=1,\ldots,r$. Note that since $s_0$ has expected order one, all $\alpha_j$ differ from zero. Now we can formulate the contribution of $E_i$ to the residue of \ZtopIg\ at $s_0$, as
\begin{equation*}
\begin{cases}
\ds\frac{1}{N_i}\left(2-m+\sum\nolimits_{j=1}^r\frac{m_j}{\alpha_j}\right),&\text{if $E_i$ is an exceptional divisor, and}\\[+3ex]
\ds\frac{1}{N_i}\sum_{j=1}^r\frac{m_j}{\alpha_j},&\text{otherwise.}
\end{cases}
\end{equation*}

It follows from the observations above that in concrete examples (in dimension two), the intersection diagram of an embedded resolution\footnote{This is the intersection diagram of the divisor $F+G$ that arises from the resolution (see Definition~\ref{def_er}).} of \Ig, together with the numerical data associated to the prime divisors, allows us to decide quickly which candidate poles are poles of \ZtopIg, and which are not.

\subsection{A formula for the topological zeta function of a monomial ideal and a non-degenerated polynomial volume form}\label{aformula}
%
%
%
We introduce a combinatorial formula for the topological zeta function of a monomial ideal and a volume form $gdx$, where $g$ is a polynomial that is non-degenerated over \C\ with respect to its global Newton polyhedron. A simplified version of the formula will be used in Section~\ref{positiveresults}.

In \cite{HMY07}, Howald, Musta{\c{t}}{\u{a}}, and Yuen provide a formula in terms of Newton polyhedra for Igusa's $p$-adic zeta function associated to a monomial ideal $\I\lhd\Z[x_1,\ldots,x_n]$. Letting it meet the older and similar one of Denef and Hoornaert \cite{DH01} for Igusa's local zeta function of a single non-degenerated polynomial, in \cite[Theorem~2.5]{BorIZFs}, we obtained a formula for the $p$-adic zeta function \ZIg\ of a monomial ideal $\I\lhd\Zp[x_1,\ldots,x_n]$, where \Zp\ denotes the ring of $p$-adic integers, and a volume form $gdx$ on $\A^n(\Qp)$, where $g$ is a polynomial over \Zp\ that is non-degenerated over \Fp\ with respect to its Newton polyhedron. Considering the motivic analogue of this formula, and specializing it to the topological version, one proves a similar formula for the topological zeta function \ZtopIg, in the case that \I\ is monomial and $g$ is non-degenerated over \C\ with respect to its global Newton polyhedron. We state this formula in Theorem~\ref{formZtopIgNP} below.

In fact the formula is heuristically obtained from the one in \cite[Theorem~2.5]{BorIZFs} by letting $p$ tend to one. The formula is also the analogue of a formula Denef and Loeser obtained in \cite[Th\'{e}or\`{e}me~5.3]{DL92} for the topological zeta function of a single non-degenerated polynomial. Their direct proof, associating an embedded resolution of singularities to the polynomial's Newton polyhedron, could also be adapted to demonstrate Theorem~\ref{formZtopIgNP}.

In order to state the theorem, we recall some definitions and notations from \cite{DH01,HMY07,BorIZFs,DL92}. For $\omega=(\omega_1,\ldots,\omega_n)\in\N^n$, we denote by $x^{\omega}$ the corresponding monomial $x_1^{\omega_1}\cdots x_n^{\omega_n}$ in $\C[x]=\C[x_1,\ldots,x_n]$. Let \I\ be a nonzero proper monomial ideal of $\C[x]$ and $g=\sum_{\omega}a_{\omega}x^{\omega}$ a nonzero polynomial in $\C[x]$, satisfying $g(0)=0$. Put $\Rplus=\{x\in\R\mid x\geqslant0\}$ and denote by $\supp(g)=\{\omega\in\N^n\mid a_{\omega}\neq0\}$ the support of $g$. The global Newton polyhedron $\Gglg$ of $g$ is defined as the convex hull of $\supp(g)$ in \Rplusn; the convex set $\Gg=\Gglg+\Rplusn$ is called the (local) Newton polyhedron of $g$. Similarly we define the Newton polyhedron \GI\ of \I\ as the convex hull in \Rplusn\ of the set $\{\omega\in\N^n\mid x^{\omega}\in \I\}$.

We say that $g$ is non-degenerated over \C\ with respect to (all the faces of) its global Newton polyhedron \Gglg, if for every face\footnote{By a face of a polyhedron we mean the polyhedron itself or one of its proper faces, which are the intersections of the polyhedron with a supporting hyperplane. See, e.g., \cite{Roc70}.} $\tau$ of \Gglg, the zero locus of $g_{\tau}=\sum_{\omega\in\tau}a_{\omega}x^{\omega}$ has no singularities in \Ccrossn, i.e., if $g_{\tau}$ and its partial derivatives $\partial g_{\tau}/\partial x_i$; $i=1,\ldots,n$; have no common root in \Ccrossn. We say that $g$ is non-degenerated over \C\ with respect to the compact faces of its (local) Newton polyhedron \Gg, if for every compact face $\tau$ of \Gg, the zero locus of $g_{\tau}$ has no singularities in \Ccrossn.

For $k_1,\ldots,k_r\in\Rn\setminus\{0\}$, we call $\Delta=\{\lambda_1k_1+\lambda_2k_2+\cdots+\lambda_rk_r\mid \lambda_i\in\Rplusnul\}\subset\Rn$ the (convex) cone strictly positively spanned by the vectors $k_1,\ldots,k_r$. When the $k_1,\ldots,k_r$ are in \Zn, we call it a rational cone. If $k_1,\ldots,k_r$ are linearly independent over \R, then $\Delta$ is called a simplicial cone. If $\Delta$ is rational and $\{k_1,\ldots,k_r\}$ is a subset of a \Z-module basis of \Zn, we call $\Delta$ a simple cone. There always exists a finite partition of $\Delta$ into cones $\Delta_i$, such that each $\Delta_i$ is strictly positively spanned by a \R-linearly independent subset of $\{k_1,\ldots,k_r\}$. We call such a decomposition a simplicial decomposition of $\Delta$ without introducing new rays. If $\Delta$ is a rational simplicial cone, we can partition $\Delta$ into a finite number of simple cones. (In general, such a decomposition requires the introduction of new rays.)

For $k\in\Rplusn$, put $\mI(k)=\inf\{k\cdot x\mid x\in\GI\}$, where $\cdot$ denotes the usual inner product on $\R^n$. We define the first meet locus $F_{\I}(k)$ of $k$ as the face $F_{\I}(k)=\{x\in\GI\mid k\cdot x=\mI(k)\}$ of \GI. For a face $\tau$ of \GI, we call $\DIt=\{k\in\Rplusn\mid F_{\I}(k)=\tau\}$ the cone associated to $\tau$. If $\tau=\GI$, we have $\DIt=\{0\}$; otherwise let $\gamma_1,\ldots,\gamma_r$ be the facets of \GI\ that contain $\tau$, and denote by $k_1,\ldots,k_r$ the unique primitive\footnote{We call a vector $k\in\Rn$ primitive if its components are integers with greatest common divisor one.} vectors in $\N^n\setminus\{0\}$ perpendicular to $\gamma_1,\ldots,\gamma_r$, respectively, then \DIt\ is the cone strictly positively spanned by $k_1,\ldots,k_r$. Note that \DIt\ has dimension $n-\dim\tau$, and that the function \mI\ is linear on the cone's topological closure. The cones \DIt, being the classes of the equivalence relation on \Rplusn\ that considers $k$ and $k'$ equivalent if $F_{\I}(k)=F_{\I}(k')$, form a partition of \Rplusn\ that we denote by \DI. We have of course analogous definitions and results for $m_g(k)$, $F_g(k)$, and \Dg, associated to $g$.

To both \I\ and $g$, we associate a partition \DIg\ of \Rplusn, consisting of all the non-empty intersections of cones in \DI\ with cones in \Dg; i.e.,
\begin{multline*}
\DIg=\{\DIt\cap\Dg(\tau')\mid\\*
\tau\textrm{ is a face of }\GI,\ \tau'\textrm{ is a face of }\Gg\textrm{ and }\DIt\cap\Dg(\tau')\neq\emptyset\}.
\end{multline*}
Equivalently, \DIg\ can be obtained as the quotient of \Rplusn\ by the equivalence relation considering $k$ and $k'$ equivalent if and only if $F_{\I}(k)=F_{\I}(k')$ and $F_g(k)=F_g(k')$. For a cone $\delta\in\DIg$, we denote by $\tau_{\delta}$ the unique face of \Gg, such that $\delta$ can be written as $\delta=\DIt\cap\Dg(\tau_{\delta})$ for some face $\tau$ of \GI.

For $k_1,\ldots,k_r$ \Q-linearly independent vectors in \Zn, we define the multiplicity of $k_1,\ldots,k_r$ as the index of the lattice $\Z k_1+\cdots+\Z k_r$ in the group of points with integral coordinates in the subspace spanned by $k_1,\ldots,k_r$ of the \Q-vector space $\Q^n$. This number $\mult(k_1,\ldots,k_r)$ equals the cardinality of the set $\Zn\cap\{\sum_{i=1}^r\lambda_ik_i\mid 0\leqslant\lambda_i<1\text{ for }i=1,\ldots,r\}$. Note that a rational simplicial cone is simple if and only if the multiplicity of its primitive generators equals one. Let $\tau$ be a face of \Gg. We define the volume $\Vol(\tau)$ of $\tau$ as follows. If $\dim\tau=0$, we put $\Vol(\tau)=1$. Otherwise, we define $\Vol(\tau)$ as the volume of $\tau\cap\Gglg=\Gamma^{\mathrm{gl}}_{g_{\tau}}$, with respect to the volume form on the affine hull\footnote{The affine hull of a subset $S\subset\Rn$ is defined as the intersection of all affine subsets of \Rn\ that contain $S$.} $\aff(\tau)$ of $\tau$, normalized such that the parallelepiped, spanned by a lattice basis of $\Zn\cap\aff(\tau)$ has volume one.

Let $k_1,\ldots,k_r$ be primitive vectors in \Nn, and let $\delta$ be the rational cone strictly positively spanned by them. We associate in the following way, to \I, $g$, and $\delta$, a rational function $J_{\I,g,\delta}$ with integer coefficients in one complex variable. If $k_1,\ldots,k_r$ are linearly independent over \R, the function is defined as
\begin{equation*}
J_{\I,g,\delta}:s\mapsto J_{\I,g,\delta}(s)=\frac{\mult(k_1,\ldots,k_r)}{\prod_{j=1}^r(\mI(k_j)s+m_g(k_j)+\sigma(k_j))},
\end{equation*}
where $\sigma(k_j)$ denotes the sum of $k_j$'s components. In the other case, we consider a partition of $\delta$ into rational simplicial cones $\delta_i$, $i\in I$, without introducing new rays, and define $J_{\I,g,\delta}$ as
\begin{equation*}
J_{\I,g,\delta}=\sum_{\substack{i\in I\\\dim\delta_i=\dim\delta}}J_{\I,g,\delta_i}.
\end{equation*}
In \cite[Lemme~5.1.1]{DL92}, one proves that the definition of $J_{\I,g,\delta}$ does not depend on the chosen partition. Finally, note that $J_{\I,g,\delta}$ is identically one for $\delta=\{(0,\ldots,0)\}$.

\begin{theorem}\label{formZtopIgNP}
Let \I\ be a nonzero proper monomial ideal of $\C[x]=\C[x_1,\ldots,x_n]$ and $g$ a nonzero polynomial in $\C[x]$, satisfying $g(0)=0$. If $g$ is non-degenerated over \C\ with respect to its global Newton polyhedron, then the global topological zeta function of \I\ and $gdx$, is given by
\begin{multline*}
\ZtopglIgs=\sum_{\substack{\delta\in\DIg\\\dim\delta=n}}J_{\I,g,\delta}(s)\\+\frac{1}{2}\sum_{\substack{\delta\in\DIg\\\dim\delta=\dim\Dg(\tau_{\delta})<n}}(-1)^{\codim\delta}(\codim\delta)!\Vol(\tau_{\delta})J_{\I,g,\delta}(s).
\end{multline*}
Note that $\codim\delta=n-\dim\delta=\dim\tau_{\delta}$.

For $g$ non-degenerated over \C\ with respect to the compact faces of its local Newton polyhedron \Gg, we obtain a local analogue of this formula by summing only over the $\delta\in\DIg$ that are not contained in any coordinate hyperplane, these are precisely the cones associated to compact faces of \GI\ and \Gg.
\end{theorem}

\begin{remark}\label{simpel}
For $g=a_{\omega}x^{\omega}$ a monomial (a special case that we will need in Section~\ref{positiveresults}), the above formula simplifies drastically to
\begin{equation*}
\ZtopglIgs=\ZtopoIgs=\sum_{\substack{\delta\in\DI\\\dim\delta=n}}J_{\I,g,\delta}(s),
\end{equation*}
where $\ZtopoIg=\ZtopIg$ denotes the (local) topological zeta function. Moreover, $m_g(k)=\omega\cdot k$ for all $k\in\Nn$ in this case.
\end{remark}

\begin{remark}\label{conj2monidtop}
We mentioned in the introduction that Howald, Musta{\c{t}}{\u{a}}, and Yuen proved the conjecture of Igusa--Denef--Loeser for Igusa's zeta function of a monomial ideal; we explain why this is not different for the topological zeta function. Considering the formula in the previous remark for $g=1$, we find that every pole of \ZtopI\ has the form $-\sigma(k)/m_I(k)$ for some primitive generator $k$ of a ray in \DI. In \cite{HMY07}, one explains that for every such ratio, there is a torus-invariant divisor $E$ on the normalized blowing-up of $\A^n$ along \I, with numerical data $(N,\nu)$ with respect to $(\I,dx)$, such that $\sigma(k)/m_I(k)=\nu/N$. From the description of the roots of the Bernstein--Sato polynomial of a monomial ideal, provided in \cite{BMS06_2}, it then follows that every pole of \ZtopI\ is a root of \bI, confirming the conjecture.
\end{remark}

\subsection{Bernstein--Sato polynomial of an ideal}
\begin{definition}[Bernstein--Sato polynomial]
Let $f$ be a polynomial in $\C[x]=\C[x_1,\ldots,x_n]$. Denote by
$\mathcal{D}=\mathcal{D}_{\A^n(\C)}$ the ring of algebraic differential operators on $\A^n(\C)$:
\begin{equation*}
\mathcal{D}=\C[x][\partial/\partial x_1,\ldots,\partial/\partial x_n].
\end{equation*}
Let $s$ be a formal indeterminate. The group
\begin{equation*}
\C[x][1/f][s]f^s=\{\varphi f^s\mid\varphi\in\C[x][1/f][s]\}
\end{equation*}
can be given the structure of a $\mathcal{D}[s]$-module by putting $(\partial/\partial x_i)f^s=s(\partial
f/\partial x_i)f^{s-1}$ for $i=1,\ldots,n$. The Bernstein--Sato polynomial (or $b$-function) $b_f(s)$ of $f$ is defined as the
monic polynomial of minimal degree in $\C[s]$, for which there exists an operator $P\in\mathcal{D}[s]$, such that
\begin{equation*}
Pf^{s+1}=b_f(s)f^s
\end{equation*}
in $\C[x][1/f][s]f^s$.

In \cite{BMS06bis}, one generalizes this notion to an arbitrary ideal
$\I=(f_1,\ldots,f_r)$ of $\C[x]$. Let $s_1,\ldots,s_r$ be formal indeterminates. The group
\begin{equation*}
M=\C[x][1/f_1\cdots f_r][s_1,\ldots,s_r]\prod_{j=1}^rf_j^{s_j}
\end{equation*}
can be turned into a $\mathcal{D}$-module by letting $(\partial/\partial x_i)$ act on $\prod_jf_j^{s_j}$ in the expected way. For $k=1,\ldots,r$, define a $\mathcal{D}$-linear action of $t_k$ on $M$ by putting $t_k(s_j)=s_j+1$, if $j=k$, and $t_k(s_j)=s_j$, otherwise, for $j=1,\ldots,r$. More precisely, for any polynomial $m$ in $r$ variables over $\C[x][1/f_1\cdots f_r]$, we have that
\begin{equation*}
t_km(s_1,\ldots,s_r)\prod_{j=1}^rf_j^{s_j}=m(s_1,\ldots,s_{k-1},s_k+1,s_{k+1},\ldots,s_r)f_k\prod_{j=1}^rf_j^{s_j}.
\end{equation*}
Note that the action of $t_k$ is bijective. Put $s=\sum_js_j$ and $s_{j,k}=s_jt_j^{-1}t_k$ for $j,k=1,\ldots,r$, and denote by
$\mathcal{D}[s_{j,k}]_{j,k}$ the ring generated by $\mathcal{D}$ and the $s_{j,k}$. Then $\mathcal{D}[s_{j,k}]_{j,k}$ acts naturally on $M$. The Bernstein--Sato polynomial \bIs\ of \I\ is defined as the monic polynomial of minimal degree in $\C[s]$, for which there exist operators $P_1,\ldots,P_r\in\mathcal{D}[s_{j,k}]_{j,k}$, such that
\begin{equation*}
\sum_{k=1}^rP_kt_k\prod_{j=1}^rf_j^{s_j}=\bIs\prod_{j=1}^rf_j^{s_j}
\end{equation*}
in $\C[x][1/f_1\cdots f_r][s_1,\ldots,s_r]\prod_{j=1}^rf_j^{s_j}$.
\end{definition}

\begin{remark}
In the case of a single polynomial, the existence and uniqueness
of the Bernstein--Sato polynomial was proved by Bernstein \cite{Ber72}
and Sato \cite{SS72}. In the case of an ideal, this follows from
the theory of the $V$-filtration of Kashiwara and Malgrange.
Budur, Musta{\c{t}}{\v{a}}, and Saito \cite{BMS06bis} proved that
the definition of \bI\ does not depend on the choice of
generators for the ideal \I. We refer to \cite{Sai06} for an
introduction to the theory of $b$-functions.
\end{remark}

The paper \cite{BMS06} by Budur, Musta{\c{t}}{\v{a}}, and Saito
provides a nice combinatorial description of the roots of the
Bernstein--Sato polynomial in the case of a monomial ideal. We will use this description in terms of Newton
polyhedra several times in this text.

\subsection{Formulation of the problem}
The following question is the subject of this paper. Although we will be dealing mostly with the dimension two case, we formulate the question for arbitrary dimension.

\begin{question}\label{vraag}
Let \I\ be a monomial ideal of $\C[x]=\C[x_1,\ldots,x_n]$ and $\bI$ its
associated Bernstein--Sato polynomial. Is it possible to find a
collection $\mathcal{G}$ of polynomials $g\in\C[x]$, such
that, for all $g\in\mathcal{G}$, every pole of the topological
zeta function \ZtopIg\ is a root of $\bI$, and such that every
root of $\bI$ is a pole of \ZtopIg\ for some $g\in\mathcal{G}$?
In other words, does $\mathcal{G}\subset\C[x]$ exist, such
that
\begin{equation}\label{eqofsets}
\bigcup_{g\in\mathcal{G}}\{s\in\C\mid s\text{ is a pole of
}\ZtopIg\}=\bI^{-1}(0)\text{?}
\end{equation}
\end{question}

As an example, we show that for principal monomial ideals \I, one obtains a positive answer to the question above. Choose $\omega=(\omega_1,\ldots,\omega_n)\in\Nn$ and let $\I=(x^{\omega})$ be the principal ideal of $\C[x_1,\ldots,x_n]$ generated by the monomial $x^{\omega}=x_1^{\omega_1}\cdots x_n^{\omega_n}$. It is well-known that the Bernstein--Sato polynomial \bI\ of $x^{\omega}$ is given by
\begin{equation*}
\bIs=\prod_{i=1}^n\prod_{j=1}^{\omega_i}\left(s+\frac{j}{\omega_i}\right).
\end{equation*}

From the formula obtained in Remark~\ref{simpel}, it follows that the topological zeta function \ZtopIg, associated to $\I=(x^{\omega})$ and a monomial $g=x^{\gamma}$, is explicitly given by
\begin{equation*}
\ZtopIgs=\frac{1}{(\omega_1s+\gamma_1+1)\cdots(\omega_ns+\gamma_n+1)},
\end{equation*}
while its poles are given by $-(\gamma_i+1)/\omega_i$; $i=1,\ldots,n$. Now we solve the question easily, considering---for example---the following collection of monomials:
\begin{equation*}
\mathcal{G}=\{1,x_1,x_1^2,\ldots,x_1^{\omega_1-1}\}\cup\{x_2,x_2^2,\ldots,x_2^{\omega_2-2}\}\cup\cdots\cup\{x_n,x_n^2,\ldots,x_n^{\omega_n-2}\}.
\end{equation*}

In the general case, we will give a negative answer to Question~\ref{vraag}. More
precisely, we will show the following facts.
\begin{enumerate}
\item In general, it is not possible to attain every root of \bI\ as a pole of \ZtopIg\ for some $g\in\C[x]$. This means that the inclusion \lq$\supset$\rq\ of Equality~\eqref{eqofsets} does not hold in general. We present a counterexample in dimension two, where \I\ is generated by three monomials.
\item In dimension two, when we restrict to ideals \I\ that can be generated by one or two monomials, every root of \bI\ can be attained as a pole of \ZtopIg\ for some $g\in\C[x]$. This is even possible when we only consider monomials $g$.
\item Attaining all the roots of \bI\ in (ii), can generally not be done without creating \lq unwanted\rq\ poles. This is, there exist ideals \I, generated by two monomials in two variables, and roots $s$ of \bI, such that every polynomial $g\in\C[x_1,x_2]$ that makes $s$ a pole of \ZtopIg, gives rise to \lq bad\rq\ poles, i.e., poles of \ZtopIg\ that are not roots of \bI. Thus for ideals in dimension two, generated by two monomials, Inclusion~\lq$\supset$\rq\ of \eqref{eqofsets} holds, but then \lq$\subset$\rq\ fails in general.
\end{enumerate}
In Sections~\ref{sect2} and \ref{sect4}, we provide counterexamples illustrating (i) and (iii), respectively. The positive result (ii) will be discussed in Section~\ref{positiveresults}.

\section{Generally, not all roots of \bI\ can be poles of \ZtopIg}\label{sect2}
To prove the assertion in the title, we consider the ideal $\I=(xy^5,x^3y^2,x^4y)\lhd\C[x,y]$, which is also mentioned in the paper of Budur, Musta{\c{t}}{\v{a}}, and Saito \cite[Example~4.3]{BMS06}. Using the authors' combinatorial description for monomial ideals, we find easily that the roots of the Bernstein--Sato polynomial \bI\ of \I, are given by
\begin{equation*}
-\frac{i}{13};\ i=5,\ldots,17;\qquad\text{and}\qquad-\frac{j}{5};\ j=2,\ldots,6.
\end{equation*}
We will focus on the root $-6/13$ and prove that there exists no $g\in\C[x,y]$, such that $-6/13$ is a pole of \ZtopIg.

To know the poles of \ZtopIg\ for a general $g$, we need an embedded resolution of \Ig\ in some open neighborhood of the origin in $\A^2(\C)$. We obtain such a local resolution by first principalizing \I\ (the origin is the only point of $\A^2(\C)$ where \I\ is not yet locally principal), and next resolving singularities and non-normal crossings of the strict transform of $g^{-1}(0)$ in the exceptional locus of the principalization.

\begin{figure}
\centering
\psset{xunit=.04545\textwidth,yunit=.04167\textwidth}
\subfigure[Numerical data for $g=1$, and, more generally, for $g(0)\neq0$.]{
\begin{pspicture}(10,6)
{
\psset{linewidth=3\pslinewidth,arrows=c-c}
\psline(0,1)(6,1)\rput[t](3,.8){$E_1(5,2)$}
\psline(5,0)(5,6)\rput[b]{90}(4.8,3){$E_3(13,5)$}
\psline(4,5)(10,5)\rput[b](7,5.2){$E_2(7,3)$}
}
\psline(1,6)(1,0)\rput[br]{90}(.8,6){$E(1,1)$}
\psline(9,6)(9,0)\rput[tl]{90}(9.2,0){$E'(1,1)$}
\end{pspicture}
}
\hfill
\subfigure[Fixed $N,N',N_i$ and minimal $\nu,\nu',\nu_i$ in case $g(0)=0$.]{
\begin{pspicture}(10,6)
{
\psset{linewidth=3\pslinewidth,arrows=c-c}
\psline(0,1)(6,1)\rput[t](3,.8){$E_1(5,3)$}
\psline(5,0)(5,6)\rput[b]{90}(4.8,3){$E_3(13,7)$}
\psline(4,5)(10,5)\rput[b](7,5.2){$E_2(7,4)$}
}
\psline(1,6)(1,0)\rput[br]{90}(.8,6){$E(1,1)$}
\psline(9,6)(9,0)\rput[tl]{90}(9.2,0){$E'(1,1)$}
\end{pspicture}
}
\caption{Intersection diagram of the minimal principalization of $\I=(xy)(y^4,x^2y,x^3)$.}\label{figprincIvbSect1samen}
\end{figure}
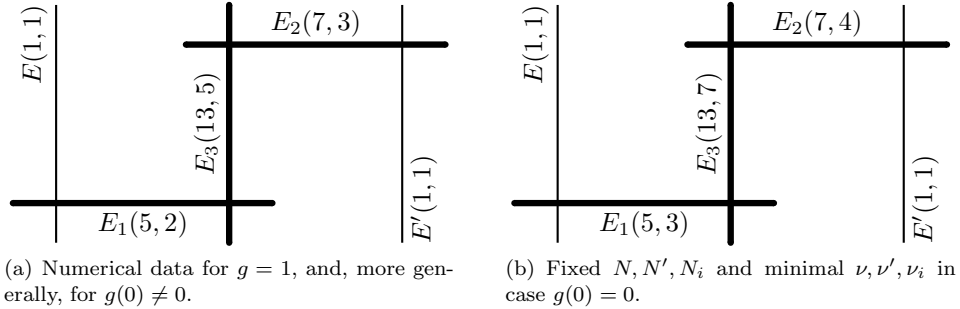

Figure~\ref{figprincIvbSect1samen}(a) shows the intersection diagram of the minimal principalization of $\I=(xy)(y^4,x^2y,x^3)$, which is a composition $h=h_1\circ h_2\circ h_3:Y=Y_3\to Y_0=\A^2(\C)$ of three blowing-ups $h_i:Y_i\to Y_{i-1}$ in points $P_{i-1}\in Y_{i-1}$. The numerical data shown are in the case $g=1$. Divisors $E$ and $E'$ denote the strict transforms of $\{y=0\}$ and $\{x=0\}$, respectively, in all $Y_i$, and $E_i$ denotes the exceptional divisor in $Y_i$ of the $i$th blowing-up $h_i$, as well as its strict transforms in $Y_j$, $j>i$, for $i=1,2,3$.

%
%
%
%

Let us denote the numerical data of $E,E',E_i$ by $(N,\nu),(N',\nu'),(N_i,\nu_i)$, respectively. While the $\nu,\nu',\nu_i$ depend on $g$, the $N,N',N_i$ only depend on \I, and are, independently of $g$, as shown in Figure~\ref{figprincIvbSect1samen}(a). If the zero locus of $g$ does not contain the origin, then the minimal principalization of \I\ is already a local embedded resolution of \Ig, and all numerical data are as shown in Figure~\ref{figprincIvbSect1samen}(a). The number $-6/13$ is not a candidate pole, and so we're done.

Suppose from now on that $g(0)=0$ and write $g=x^ny^m\tilde{g}$, with $\tilde{g}$ a polynomial not divisible by $x$ or $y$. Then $\nu=1+m\geqslant1$, $\nu'=1+n\geqslant1$, and $\nu_1=\nu+\nu'+\mult_{P_0}(\tilde{g})\geqslant3$, with $P_0=(0,0)\in\A^2(\C)$, the center of the first blowing-up. It follows that $\nu_2=\nu'+\nu_1+\mult_{P_1}(\tilde{g})\geqslant4$ and $\nu_3=\nu_1+\nu_2+\mult_{P_2}(\tilde{g})\geqslant7$, where $P_1$ is the intersection point of $E'$ and $E_1$ in $Y_1$, $P_2$ is the intersection point of $E_1$ and $E_2$ in $Y_2$, and $\tilde{g}$ denotes---by abuse of notation---also the strict transforms of $\tilde{g}$ in $\mathcal{O}_{Y_1}(Y_1)$ and $\mathcal{O}_{Y_2}(Y_2)$. The minimal principalization diagram with minimal $\nu,\nu',\nu_i$, is shown in Figure~\ref{figprincIvbSect1samen}(b).

Notice that all candidate poles, arising from the minimal principalization, are strictly smaller than $-6/13$. However, we may not yet have a local embedded resolution of \Ig. Generally, we still need to perform a series of blowing-ups with centers in the exceptional locus. Two types of blowing-ups must be considered. If we blow up in a point $P$ on an exceptional divisor $E_i(N_i,\nu_i)$, that is not contained in any other relevant divisor (exceptional, $E$ or $E'$), then the resulting exceptional divisor $E_j$ will have numerical data $(N_i,\nu_i+1+\mult_P(\tilde{g}))$, and the new candidate pole will be strictly smaller than the one belonging to $E_i$. If we blow up in the intersection point $P$ of an exceptional divisor $E_i(N_i,\nu_i)$ and another relevant divisor $E_{\ast}(N_{\ast},\nu_{\ast})$, then the resulting exceptional divisor $E_j$ will have numerical data $(N_i+N_{\ast},\nu_i+\nu_{\ast}+\mult_P(\tilde{g}))$, and the new candidate pole will be strictly\footnote{Although this is not necessary for the current argumentation, we may assume this inequality to be strict, because, if $\mult_P(\tilde{g})$ were zero, there would be no reason for this blowing-up.} smaller than the maximum of the candidate poles belonging to $E_i$ and $E_{\ast}$.

Those facts imply that the maximum of our set of candidate poles will not increase by performing the blowing-ups that are still needed to obtain a local embedded resolution of \Ig, and so none of the candidate poles will equal $-6/13$.

\begin{remark}Using {\sc Singular}'s {\tt gaussman.lib} library \cite{sing,singlibbernstein}, we verified that $-6/13$ is also a root of the Bernstein--Sato polynomial $b_f$ of $f=xy^5+x^3y^2+x^4y\in\I$. The series of blowing-ups that principalized \I, also gives an embedded resolution of $f^{-1}(0)\subset\A^n(\C)$, having a single singularity in the origin. (We call $f^{-1}(0)$ a generic curve of \I.) Therefore, the same argumentation yields a counterexample for the analogue of Question~\ref{vraag} in the case of one polynomial as well.
\end{remark}

\begin{remark}
It may be an interesting question for further study to ask, given a monomial (or arbitrary) ideal \I, which roots of \bI\ can be obtained as a pole of \ZtopIg\ for some $g$, and which cannot, and what is the \lq nature\rq\ of the latter?
\end{remark}

\section{A partially positive result for ideals in dimension two generated by two monomials}\label{positiveresults}
For ideals \I\ generated by two monomials in two variables, we show that every root of \bI\ can be obtained as a pole of \ZtopIg\ for some polynomial $g$. This proves one inclusion of \eqref{eqofsets} from Question~\ref{vraag} in the present case. In dimension two, this result is in some sense optimal: on the one hand it follows from Section~\ref{sect2} that the result cannot be generalized to more than two generators; on the other hand we will see in Section~\ref{sect4} that we cannot expect the other inclusion of \eqref{eqofsets} to hold in general (in the present case).

%
%
\begin{figure}
\centering
\psset{xunit=.0284\textwidth,yunit=.0284\textwidth}
\subfigure[Newton polyhedron $\GI$.]{
\begin{pspicture}(-1,-1.4)(10.5,10.5)
\pspolygon*[linecolor=lightgray,linewidth=0pt](2,10.5)(2,5)(6,3)(10.5,3)(10.5,10.5)
\psaxes[labels=none,ticks=none]{->}(0,0)(0,0)(10.5,10.5)
\psdots(2,5)(6,3) \psline(2,10.5)(2,5)(6,3)(10.5,3)
\psline[linestyle=dotted](0,5)(2,5)(2,0)
\psline[linestyle=dotted](0,3)(6,3)(6,0)
\uput[l](2,7.75){$\tau_1$} \uput[230](2,5){$\tau_2$}
\uput[240](4,4){$\tau_3$} \uput[244](6,3){$\tau_4$}
\uput[d](8.25,3){$\tau_5$} \uput[40](2,5){$(a,b)$}
\uput[66](6,3){$(c,d)$} \rput(6.5,7){$\GI=\tau_0$}
\uput[l](0,5){$b$}\uput[l](0,3){$d$}
\uput[d](2,0){$\phantom{-}a\phantom{-}$}\uput[d](6,0){$\phantom{-}c\phantom{-}$}
\uput[dl](0,0){$\phantom{-}0$}
\end{pspicture}
}
\hfill
\subfigure[Dimension $2$ cones in $\DI$.]{
\begin{pspicture}(-3.4,-1.4)(8.5,10.5)
\psaxes[labels=none,ticks=none](0,0)(0,0)(8.5,10.5)
\psline(0,0)(5.25,10.5) \psline[linewidth=2pt]{C->}(0,0.03)(0,2.5)
\psline[linewidth=2pt]{->}(0,0)(3,6)
\psline[linewidth=2pt]{C->}(0.03,0)(2.5,0)
\rput(6.5,5.67){$\DI(\tau_2)$} \rput(1.75,7){$\DI(\tau_4)$}
\uput[l](0,1.25){$v_3(0,1)$}
\uput[330](1.5,3){\parbox{.185\textwidth}{$v_2\left(\frac{b-d}{e},\frac{c-a}{e}\right),$\\{\scriptsize
$e=\gcd(b-d,c-a)$}}} \uput[d](1.25,0){$v_1(1,0)$}
\end{pspicture}
}
\hfill
\subfigure[$S_{\tau_3}$.]{
\begin{pspicture}(-2.5,-1.4)(6.7,10.5)
\pspolygon*[linecolor=lightgray,linewidth=0pt](0,0)(4,0)(6,3)(6,5)(2,5)(0,2)
\psaxes[labels=none,ticks=none]{->}(0,0)(0,0)(6.7,10.5)
\psline[linestyle=dotted](0,5)(2,5)(2,0)
\psline[linestyle=dotted](0,3)(6,3)(6,0)
\psline[linewidth=1.4pt](0,2)(0,0)(4,0)(6,3)
\psline[linewidth=.2pt,doubleline=true,doublesep=1pt,fillcolor=white](6,3)(6,5)(2,5)(0,2)
\psdots[dotstyle=o,dotsize=2.5pt](0,2)(6,3)
\psdots[dotsize=2.5pt](0,0)(0,1)(1,0)(1,1)(1,2)(1,3)(2,0)(2,1)(2,2)(2,3)(2,4)(3,0)(3,1)(3,2)(3,3)(3,4)(4,0)(4,1)(4,2)(4,3)(4,4)(5,2)(5,3)(5,4)
\uput[l](0,5){$b$}\uput[l](0,3){$d$}\uput[l](0,2){$b-d$}
\uput[d](2,0){$\phantom{-}a\phantom{-}$}\uput[d](6,0){$\phantom{-}c\phantom{-}$}\uput[d](4,0){$c-a$}
\uput[dl](0,0){$\phantom{-}0$}
\psframe*[linecolor=white](2.3,1.3)(3.7,2.7)
\rput(3,2){$S_{\tau_3}$}
\end{pspicture}
}
\caption{$\GI$, $\DI$, and $S_{\tau_3}$, associated to $\I=(x^ay^b,x^cy^d)\lhd\C[x,y]$.}\label{fig_samen_Sect2}
\end{figure}
%
%

Let $\I=(x^ay^b,x^cy^d)\lhd\C[x,y]$ be an ideal, generated by two
monomials. In view of the affirmative answer to
Question~\ref{vraag} in the principal ideal case, we may suppose
that \I\ cannot be generated by only one monomial; i.e., we may
assume $a<c$ and $b>d$. Figure~\ref{fig_samen_Sect2}(a--b) shows
the Newton polyhedron of \I\ and the associated partition of
\Rplusn. Following the combinatorial
description of Budur, Musta{\c{t}}{\v{a}}, and Saito \cite{BMS06},
we find the roots of \bI:
\begin{align*}
\text{roots associated to face $\tau_1$ (if $a\neq0$): }&-\frac{i+1}{a};\ i=0,\ldots,a-1;\\
\text{roots associated to face $\tau_3$: }&-\varphi(k+1,l+1);\ (k,l)\in\N^2\cap S_{\tau_3};\\
\text{roots associated to face $\tau_5$ (if $d\neq0$): }&-\frac{j+1}{d};\ j=0,\ldots,d-1;
\end{align*}
where
\begin{equation*}
\varphi:\R^2\to\R:(x,y)\mapsto\varphi(x,y)=\frac{v_2\cdot(x,y)}{v_2\cdot(a,b)}=\frac{(b-d)x+(c-a)y}{bc-ad}
\end{equation*}
is the unique linear function on $\R^2$ with rational coefficients that is identically one on $\tau_3$, and $S_{\tau_3}$ is the set of solutions $(x,y)\in\R^2$ to the following system of inequalities (cfr.\ Figure~\ref{fig_samen_Sect2}(c)):
\begin{equation*}
\left\{
\begin{alignedat}{2}
0&\leqslant x&&<c,\\
0&\leqslant y&&<b,\\
\frac{d(x+a-c)}{a}&\leqslant y&&<\frac{dx}{a}+b-d.
\end{alignedat}
\right.
\end{equation*}

We want to show that every root of \bI\ can be obtained as a pole of \ZtopIg\ for some monomial $g=x^ny^m\in\C[x,y]$. We calculate \ZtopIgs\ explicitly for $g=x^ny^m$, using the formula obtained in Remark~\ref{simpel}:
\begin{align*}
&\ZtopIgs\\*
&=\sum_{\substack{\delta\in\DI\\\dim\delta=2}}J_{\I,g,\delta}(s)\\*
&=J_{\I,g,\DI(\tau_2)}(s)+J_{\I,g,\DI(\tau_4)}(s)\\*
&=\frac{\mult(v_1,v_2)}{(\mI(v_1)s+m_g(v_1)+\sigma(v_1))(\mI(v_2)s+m_g(v_2)+\sigma(v_2))}\\*
&\qquad+\frac{\mult(v_2,v_3)}{(\mI(v_2)s+m_g(v_2)+\sigma(v_2))(\mI(v_3)s+m_g(v_3)+\sigma(v_3))}\\*
&=\frac{c-a}{[as+n+1][(a(b-d)+b(c-a))s+n(b-d)+m(c-a)+b-d+c-a]}\\*
&\qquad+\frac{b-d}{[(a(b-d)+b(c-a))s+n(b-d)+m(c-a)+b-d+c-a][ds+m+1]}.
\end{align*}

The candidate poles of \ZtopIg\ are
\begin{alignat*}{2}
s_1&=-\frac{n+1}{a}&\qquad&\text{(if $a\neq0$)},\\
s_2&=-\frac{(b-d)(n+1)+(c-a)(m+1)}{bc-ad},&&\text{and}\\
s_3&=-\frac{m+1}{d}&&\text{(if $d\neq0$)}.
\end{alignat*}
Suppose $a\neq0$. If $s_1=s_2$, then the limit
\begin{equation*}
\lim_{s\to s_1}(s-s_1)^2\ZtopIgs=
\begin{cases}
\ds\frac{1}{bc-ad}\left(\frac{c-a}{a}+\frac{b-d}{d}\right),&\text{if $d\neq0$ and $s_1=s_3$};\\[+2ex]
\ds\frac{c-a}{a(bc-ad)},&\text{otherwise};
\end{cases}
\end{equation*}
is a positive rational number, and consequently, $s_1$ is a pole of order two. Otherwise, the residue
\begin{equation*}
\res_{s_1}\ZtopIg=
\begin{cases}
\ds\frac{1}{a(m+1)-b(n+1)}\left(1+\frac{a(b-d)}{d(c-a)}\right),&\text{if $d\neq0$ and $s_1=s_3$};\\[+2ex]
\ds\frac{1}{a(m+1)-b(n+1)},&\text{otherwise};
\end{cases}
\end{equation*}
of \ZtopIg\ at $s_1$ is a nonzero rational number; hence $s_1$ is a pole of order one. We conclude analogously for the candidate pole $s_3$ if $d\neq0$. Similarly, $s_2$ is a double pole if $s_2$ coincides with another candidate pole ($s_1$ or $s_3$); in the other case the residue at $s_2$ equals
\begin{equation*}
\res_{s_2}\ZtopIg=-\frac{(b-d)(n+1)+(c-a)(m+1)}{[a(m+1)-b(n+1)][c(m+1)-d(n+1)]}\neq0,
\end{equation*}
and $s_2$ is a simple pole. So we find that all candidate poles are effectively poles.

Comparing the roots of \bI\ to the poles of \ZtopIg\ for $g=x^ny^m$, we see that a root of the form $-(i+1)/a$, $-\varphi(k+1,l+1)$, and $-(j+1)/d$ is a pole of \ZtopIg\ for $g=x^i$, $g=x^ky^l$, and $g=y^j$, respectively, what was to be shown.

\section{Attaining all roots of \bI\ in the two monomial case can---in general---not be done without creating \lq unwanted\rq\ poles}\label{sect4}
We give an example of an ideal $\I\lhd\C[x,y]$, generated by two monomials, and a root $s_0$ of \bI, such that every $g\in\C[x,y]$ that makes $s_0$ a pole of \ZtopIg, gives rise to other poles of \ZtopIg\ that are no roots of \bI.

Consider the ideal $\I=(xy,x^5)=(x)(y,x^4)\lhd\C[x,y]$. \ZtopI\ has a single pole $-1$, while the roots of \bI\ are $-i/5$; $i=5,\ldots,9$; again following \cite{BMS06}. The roots $-1$, $-6/5$, and $-9/5$ can be obtained as poles of \ZtopIg, without unwanted poles being created. For $-1$, this is clear, for the other two, we can take $g=x+y^4$. This is not possible for root $s_0=-7/5$, as we will show.

%
%
\begin{figure}
\centering
\psset{xunit=.04545\textwidth,yunit=.04167\textwidth}
\begin{pspicture}(10,6)
{
\psset{linewidth=3\pslinewidth,arrows=c-c}
\psline(1,6)(1,0)\rput[br]{90}(.8,6){$E_4(5,5)$}
\psline(0,1)(6,1)\rput[t](3,.8){$E_3(4,4)$}
\psline(5,0)(5,6)\rput[b]{90}(4.8,3){$E_2(3,3)$}
\psline(4,5)(10,5)\rput[b](7,5.2){$E_1(2,2)$}
}
\psline(9,6)(9,0)\rput[tl]{90}(9.2,0){$E_0(1,1)$}
\end{pspicture}
\caption{Intersection diagram of the minimal principalization of $\I=(xy,x^5)$. Numerical data shown are for $g(0)\neq0$. For general $g$, we have $\nu_k=k+1+\min\{i+kj\mid(i,j)\in\supp g\}$; $k=0,\ldots,4$.}\label{figprincI}
\end{figure}
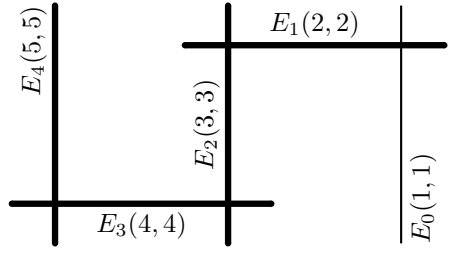
%
%

Consider a general polynomial $g=\sum_{(i,j)}c_{i,j}x^iy^j\in\C[x,y]$ and the minimal principalization of \I, whose intersection diagram and numerical data (with respect to $dx$ and $gdx$) can be found in Figure~\ref{figprincI}. We will concentrate first on those $g$ for which $\nu_4=7$ and thus \ZtopIg\ has a candidate pole $-7/5$, and prove that those $g$ cause \ZtopIg\ to have poles that are non-roots. In a second step, we will prove that \ZtopIg\ can never have $-7/5$ as a pole, if $\nu_4\neq7$, and then we're done.

\subsection*{Step I: $\nu_4=7$}
Since $\nu_4=5+\min\{i+4j\mid(i,j)\in\supp g\}$ and we may assume that $g$ has a term in $y^n$ for some $n$ (otherwise \ZtopIg\ would have a pole $\leqslant-2$), we are looking at $g$ of the form
\begin{equation*}
g=ax^2+by^n+\sum_{(i,j)\in S_n}c_{i,j}x^iy^j,
\end{equation*}
with $n\in\Nnul$; $a,b\in\C\setminus\{0\}$; and $S_n$ a finite subset of
\begin{equation*}
\N^2\setminus\{(0,0),(1,0),(2,0);\quad(0,1),(0,2),\ldots,(0,n)\}.
\end{equation*}

Now the goal is to calculate a local embedded resolution of \Ig\ (with $g$ of the above form), starting from the minimal principalization of \I, in order to determine the poles of \ZtopIg. To do so, it turns out we need to distinguish between two main cases, depending on the shape of $g$'s Newton polyhedron:
\begin{description}
\item[{\normalfont Case 1}] $c_{1,j}=0$ for all $j<n/2$;
\item[{\normalfont Case 2}] the complementary case, where
\begin{equation*}
0\neq m=\min\{j\mid(1,j)\in\supp g\}<n/2.
\end{equation*}
\end{description}
In the first case, we will show that $-(3n+2)/(n+2)$ is a pole of \ZtopIg, while in the second case, we show \ZtopIg\ has poles in $-(3m+1)/(m+1)$ and $-(2n-m+1)/(n-m+1)$. One checks that all three poles are, for all values of $n$ and $m$, outside the set of roots of \bI.

\subsubsection*{Case 1}
\begin{figure}
\centering
\psset{xunit=.03571\textwidth,yunit=.04167\textwidth}
\begin{pspicture}(28,6)
{
\psset{linewidth=3\pslinewidth,arrows=c-c}
\psline(1,6)(1,0)\rput[br]{90}(.8,6){$E_4(5,7)$}
\psline(0,1)(6,1)\rput[t](3,.8){$E_3(4,6)$}
\psline(5,0)(5,6)\rput[b]{90}(4.8,3){$E_2(3,5)$}
\psline(4,5)(10,5)\rput[b](7,5.2){$E_1(2,4)$}
\psline(9,6)(9,0)\rput[t]{90}(9.2,3){$E_1'(3,7)$}
\psline(8,1)(14,1)\rput[t](11,.8){$E_2'(4,10)$}
\psline(13,0)(13,6)\rput[b]{90}(12.8,3){$E_3'(5,13)$}
\psline(12,5)(14,5)
\psline[linestyle=dotted,dotsep=.2](14,5)(16,5)
\psline(16,5)(18,5)
\psline(17,0)(17,6)\rput[l](17.2,3){$E_{k-2}'$}
\psline(16,1)(28,1)\rput[t](22,.8){$E_{k-1}'$}
}
\psline(27,6)(27,0)\rput[tr]{90}(27.2,6){$E_0(1,1)$}
\parabola[linestyle=dashed](19,0)(22,5)\rput[b](22,5.2){$\widetilde{g=0}\ (0,2)$}
\end{pspicture}
\caption{Intersection diagram of an embedded resolution of \Ig\ for Case~1.3(a). Other numerical data are $E_i'(2+i,4+3i)$; $i=1,\ldots,k-1$.}\label{figcase13a}
\end{figure}
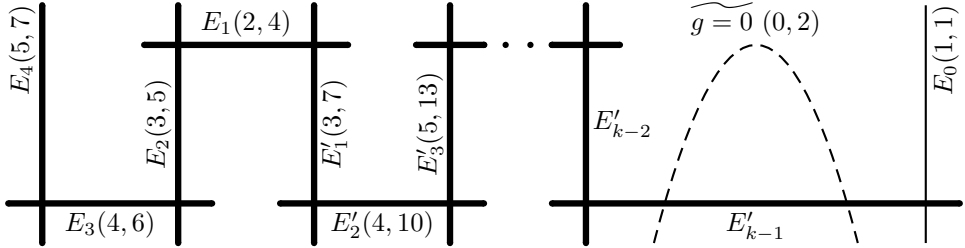

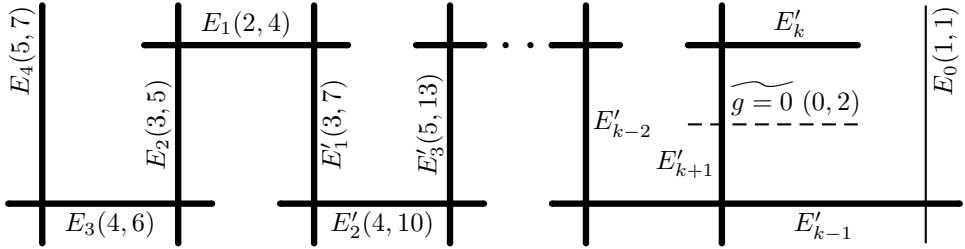
\begin{figure}
\centering
\psset{xunit=.03571\textwidth,yunit=.04167\textwidth}
\begin{pspicture}(28,6)
{
\psset{linewidth=3\pslinewidth,arrows=c-c}
\psline(1,6)(1,0)\rput[br]{90}(.8,6){$E_4(5,7)$}
\psline(0,1)(6,1)\rput[t](3,.8){$E_3(4,6)$}
\psline(5,0)(5,6)\rput[b]{90}(4.8,3){$E_2(3,5)$}
\psline(4,5)(10,5)\rput[b](7,5.2){$E_1(2,4)$}
\psline(9,6)(9,0)\rput[t]{90}(9.2,3){$E_1'(3,7)$}
\psline(8,1)(14,1)\rput[t](11,.8){$E_2'(4,10)$}
\psline(13,0)(13,6)\rput[b]{90}(12.8,3){$E_3'(5,13)$}
\psline(12,5)(14,5)
\psline[linestyle=dotted,dotsep=.2](14,5)(16,5)
\psline(16,5)(18,5)
\psline(17,0)(17,6)\rput[l](17.2,3){$E_{k-2}'$}
\psline(16,1)(28,1)\rput[t](24,.8){$E_{k-1}'$}
\psline(21,0)(21,6)\rput[r](20.8,2){$E_{k+1}'$}
\psline(20,5)(25,5)\rput[b](23,5.2){$E_k'$}
} \psline(27,6)(27,0)\rput[tr]{90}(27.2,6){$E_0(1,1)$}
\psline[linestyle=dashed](20,3)(25,3)\rput[bl](21.27,3.2){$\widetilde{g=0}\
(0,2)$}
\end{pspicture}
\caption{Intersection diagram of an embedded resolution of \Ig\ for Case~1.3(b). Other numerical data are $E_i'(2+i,4+3i)$; $i=1,\ldots,k-1$; $E_k'(k+1,3k+3)$; and $E_{k+1}'(2k+2,6k+5)$.}\label{figcase13b}
\end{figure}

Starting from the minimal principalization of \I, we perform the necessary extra blowing-ups to obtain a (local) embedded resolution of \Ig. We further distinguish three subcases:
\begin{description}
\item[{\normalfont Case 1.1}] $n=1$;
\item[{\normalfont Case 1.2}] $n=2k+1$ odd, $k\geqslant1$;
\item[{\normalfont Case 1.3}] $n=2k$ even, $k\geqslant1$.
\end{description}
In Case~1.1 the minimal principalization is already an embedded resolution. The corresponding intersection diagram and numerical data are shown in Figure~\ref{figcase11}. In Case~1.2 we need to blow up $k+1$ more times to obtain an embedded resolution. For Case~1.3 there are two possibilities: dependently, $k-1$ or $k+1$ extra blowing-ups are needed. Intersection diagrams and numerical data can be found in Figures~\ref{figcase12}--\ref{figcase13b}.

%
%
\begin{figure}
\centering
\psset{xunit=.04545\textwidth,yunit=.04167\textwidth}
\begin{pspicture}(10,6)
{
\psset{linewidth=3\pslinewidth,arrows=c-c}
\psline(1,6)(1,0)\rput[br]{90}(.8,6){$E_4(5,7)$}
\psline(0,1)(6,1)\rput[t](3,.8){$E_3(4,6)$}
\psline(5,0)(5,6)
\psline(4,5)(10,5)\rput[b](7,5.2){$E_1(2,3)$}
}
\psline(9,6)(9,0)\rput[tl]{90}(9.2,0){$E_0(1,1)$}
\psline[linestyle=dashed](2.5,3)(7.5,3)\rput[br](4.8,3.2){$\widetilde{g=0}$}\rput[bl](5.2,3.2){$(0,2)$}
\rput[l](5.2,2){$E_2(3,5)$}
\end{pspicture}
\caption{Intersection diagram of an embedded resolution of \Ig\ for Case~1.1.}\label{figcase11}
\end{figure}
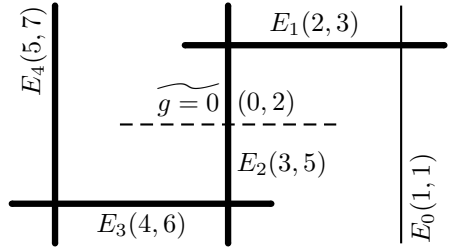
%
%

%
%
\begin{figure}
\centering
\psset{xunit=.04545\textwidth,yunit=.04167\textwidth}
\begin{pspicture}(22,6)
{
\psset{linewidth=3\pslinewidth,arrows=c-c}
\psline(1,6)(1,0)\rput[br]{90}(.8,6){$E_4(5,7)$}
\psline(0,1)(6,1)\rput[t](3,.8){$E_3(4,6)$}
\psline(5,0)(5,6)\rput[b]{90}(4.8,3){$E_2(3,5)$}
\psline(4,5)(10,5)\rput[b](7,5.2){$E_1(2,4)$}
\psline(9,6)(9,0)\rput[t]{90}(9.2,3){$E_1'(3,7)$}
\psline(8,1)(12,1)\rput[t](10.5,.8){$E_2'(4,10)$}
\psline[linestyle=dotted,dotsep=.27](12,1)(14,1)
\psline(14,1)(18,1)\rput[t](15.5,0.8){$E_{k-1}'$}
\psline(17,0)(17,6)\rput[r](16.8,2){$E_{k+1}'$}
\psline(16,5)(22,5)\rput[b](19,5.2){$E_k'$}
}
\psline(21,6)(21,0)\rput[tl]{90}(21.2,0){$E_0(1,1)$}
\psline[linestyle=dashed](14.5,3)(19.5,3)\rput[br](16.8,3.2){$\widetilde{g=0}$}\rput[bl](17.2,3.2){$(0,2)$}
\end{pspicture}
\caption{Intersection diagram of an embedded resolution of \Ig\ for Case~1.2. Other numerical data are $E_i'(2+i,4+3i)$; $i=1,\ldots,k-1$; $E_k'(k+2,3k+3)$; and $E_{k+1}'(2k+3,6k+5)$.}\label{figcase12}
\end{figure}
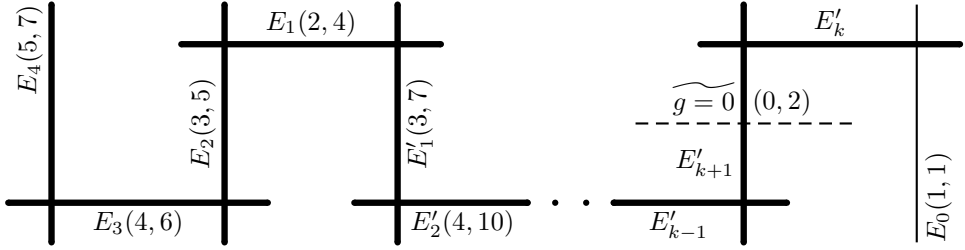
%
%

%
%
%
%

%
%
%
%

Calculating the residues of candidate poles
\begin{alignat*}{2}
-\frac{\nu_2}{N_2}&=-\frac53&\qquad&\text{(Case~1.1)},\\*
-\frac{\nu'_{k+1}}{N'_{k+1}}&=-\frac{6k+5}{2k+3}&&\text{(Case~1.2)},\quad\text{and}\\*
-\frac{\nu'_{k-1}}{N'_{k-1}}&=-\frac{3k+1}{k+1}&&\text{(Case~1.3)},
\end{alignat*}
we conclude these are the poles of \ZtopIg\ we were looking for.

\subsubsection*{Case 2}
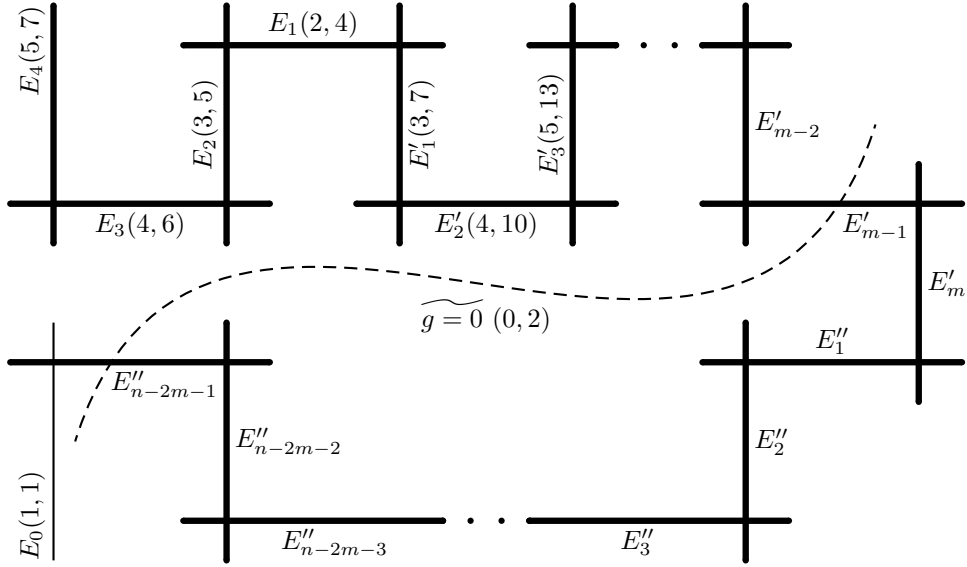
\begin{figure}
\centering
\psset{xunit=.04545\textwidth,yunit=.04167\textwidth}
\begin{pspicture}(0,-8)(22,6)
{ \psset{linewidth=3\pslinewidth,arrows=c-c}
\psline(1,6)(1,0)\rput[br]{90}(.8,6){$E_4(5,7)$}
\psline(0,1)(6,1)\rput[t](3,.8){$E_3(4,6)$}
\psline(5,0)(5,6)\rput[b]{90}(4.8,3){$E_2(3,5)$}
\psline(4,5)(10,5)\rput[b](7,5.2){$E_1(2,4)$}
\psline(9,6)(9,0)\rput[t]{90}(9.2,3){$E_1'(3,7)$}
\psline(8,1)(14,1)\rput[t](11,.8){$E_2'(4,10)$}
\psline(13,0)(13,6)\rput[b]{90}(12.8,3){$E_3'(5,13)$}
\psline(12,5)(14,5)
\psline[linestyle=dotted,dotsep=.27](14,5)(16,5)
\psline(16,5)(18,5)
\psline(17,0)(17,6)\rput[l](17.2,3){$E_{m-2}'$}
\psline(16,1)(22,1)\rput[t](20,.8){$E_{m-1}'$}
\psline(21,2)(21,-4)\rput[l](21.2,-1){$E_m'$}
\psline(22,-3)(16,-3)\rput[b](19,-2.8){$E_1''$}
\psline(17,-2)(17,-8)\rput[l](17.2,-5){$E_2''$}
\psline(18,-7)(12,-7)\rput[t](14.5,-7.2){$E_3''$}
\psline[linestyle=dotted,dotsep=.27](12,-7)(10,-7)
\psline(10,-7)(4,-7)\rput[t](7.5,-7.2){$E_{n-2m-3}''$}
\psline(5,-8)(5,-2)\rput[l](5.2,-5){$E_{n-2m-2}''$}
\psline(6,-3)(0,-3)\rput[tr](4.8,-3.2){$E_{n-2m-1}''$} }
\psline(1,-2)(1,-8)\rput[bl]{90}(.8,-8){$E_0(1,1)$}
\psbezier[linestyle=dashed](1.5,-5)(4.5,5.5)(17,-7.5)(20,3)
\rput[t](11,-1.4){$\widetilde{g=0}\ (0,2)$}
\end{pspicture}
\caption{Intersection diagram of an embedded resolution of \Ig\ for Case~2. Other numerical data are $E_i'(2+i,4+3i)$; $i=1,\ldots,m-1$; $E_m'(m+2,3m+3)$; and $E_j''(m+j+2,3m+2j+3)$; $j=1,\ldots,n-2m-1$.}\label{figcase2}
\end{figure}

After two consecutive series of (extra) blowing-ups (of lengths $m$ and $n-2m-1$), we find an embedded resolution of \Ig\ (see Figure~\ref{figcase2} for the intersection diagram and numerical data). The strict transform of $\{g=0\}$ intersects two exceptional divisors, $E'_{m-1}(m+1,3m+1)$ and $E''_{n-2m-1}(n-m+1,2n-m+1)$, resulting in two effective poles of \ZtopIg:
\begin{equation*}
-\frac{3m+1}{m+1}\qquad\text{and}\qquad-\frac{2n-m+1}{n-m+1},
\end{equation*}
the ones we announced above.

%
%
%
%

\subsection*{Step II: $\nu_4\neq7$}
Let $g$ be a polynomial such that $\nu_4\neq7$. We prove that \ZtopIg\ cannot have a pole in $-7/5$. Recall that $\nu_4=5+\min\{i+4j\mid(i,j)\in\supp g\}$. If $\nu_4=5$, then $g(0)\neq0$ and $\ZtopIg=\ZtopI$ has only one pole, namely $-1$.

If $\nu_4=6$, we know $g$ does not have a constant term and does have a term in $x$. Consequently, $\nu_0\in\{1,2\}$, $\nu_1=3$, $\nu_2=4$, and $\nu_3=5$, and all candidate poles, arising from the minimal principalization of \I, differ from $-7/5$. Moreover, the strict transform of $\{g=0\}$, under this principalization, does not intersect exceptional divisors $E_2$, $E_3$, and $E_4$, and does not intersect $E_1^{\circ}$, unless transversally (this is if and only if $g$ has a term in $y$). The strict transform may intersect $E_0^{\circ}$, but this does not affect \ZtopIg. It follows that all candidate poles, created in going---in a minimal way---from the minimal principalization to a local embedded resolution of \Ig, come from an exceptional divisor lying above $E_0\cap E_1$.

Recall from Section~\ref{sect2} that
\begin{enumerate}
\item if we blow up in a point $P$ on an exceptional divisor $E_i(N_i,\nu_i)$, that is not contained in any other relevant divisor (exceptional or $E_0$), then the resulting exceptional divisor $E_j$ will have numerical data $(N_i,\nu_i+1+\mult_P(\tilde{g}))$, and the new candidate pole will be strictly smaller than the one belonging to $E_i$, and
\item if we blow up in the intersection point $P$ of an exceptional divisor $E_i(N_i,\nu_i)$ and another relevant divisor $E_{\ast}(N_{\ast},\nu_{\ast})$, then the resulting exceptional divisor $E_j$ will have numerical data $(N_i+N_{\ast},\nu_i+\nu_{\ast}+\mult_P(\tilde{g}))$, and the new candidate pole will be strictly\footnote{Assuming $\mult_P(\tilde{g})\geqslant1$.} smaller than the maximum of the candidate poles belonging to $E_i$ and $E_{\ast}$.
\end{enumerate}
If $\nu_4\geqslant8$, it follows that $\nu_0\geqslant1$, $\nu_1\geqslant3$, $\nu_2\geqslant5$, and $\nu_3\geqslant7$, and all candidate poles, except for $-\nu_0/N_0$, arising from the minimal principalization of \I, are strictly smaller than $-7/5$. Consequently, exceptional divisors lying above $E_1^{\circ}\cup E_2\cup E_3\cup E_4$, also give candidate poles strictly smaller than $-7/5$. So, also in this case, we only have to be concerned with exceptional divisors lying above $E_0\cap E_1$.

Suppose $\nu_4\not\in\{5,7\}$. We have to check if a candidate pole, coming from an exceptional divisor lying above $E_0\cap E_1$, can equal $-7/5$. We know $-\nu_0/N_0\leqslant-1$ and $-\nu_1/N_1\leqslant-3/2<-7/5$. Blowing-up in $E_0\cap E_1$ results in a new candidate pole\footnote{We assume $\mult_{E_0\cap E_1}(\tilde{g})\geqslant1$, otherwise there would be no reason for this blowing-up.} $\leqslant-5/3<-7/5$. This makes that, from now on, we only have to consider exceptional divisors lying above the intersection of $E_0$ and this last created exceptional divisor. Blowing up in this intersection point gives a new candidate pole, being maximally $-7/4<-7/5$. Proceeding in this way, we obtain a sequence of candidate poles, the $n$th being smaller than or equal to $-(3+2n)/(2+n)$, which is strictly smaller than $-7/5$ for all $n$. We conclude $-7/5$ can never be a pole of \ZtopIg.

\bibliography{bartbib}
\bibliographystyle{amsplain}
\end{document}